\documentclass[12pt]{amsart}

\usepackage{amsmath}
\usepackage{amsfonts}
\usepackage{amssymb}
\usepackage{amsthm}
\usepackage{hyperref}
\usepackage{enumerate}
\usepackage{cleveref}
\usepackage{mathrsfs}
\usepackage[top=0.9in, bottom=1.15in, left=1.1in, right=1.1in]{geometry}
 \usepackage{hyperref}
\hypersetup{colorlinks=true,linkcolor=blue,citecolor=magenta}

\newtheorem*{remark}{Remark}

\Crefname{conjecture}{Conjecture}{Conjectures}

\theoremstyle{plain}

\theoremstyle{plain}

\author{Robert Schneider}

\address{Department of Mathematical Sciences\newline
Michigan Technological University\newline
Houghton, Michigan 49931, U.S.A.}
\email{robertsc@mtu.edu}

\title{Partition-theoretic model of prime distribution, II}
 
\begin{document}
 
\begin{abstract}
In recent work by Botkin, Dawsey, Hemmer, Just and the present author,   a {deterministic} model of prime number distribution is developed based on properties of integer partitions that gives almost exact estimates for $\pi(n)$, the number of primes less than or equal to positive integer $n$, up to $n=10{,}000$. In this follow-up paper, the author summarizes the ideas behind this partition-theoretic model of primes  and   formulates a computational model   that is  {practically exact} in its estimates of $\pi(n)$   up to $n=100,000$. 
\end{abstract}

\maketitle

\section{Introduction}\label{sect1}
\subsection{Partition model of prime distribution}\label{intro} The   sieve of Eratosthenes  provides a   method to identify prime numbers out of the larger set of positive integers $\mathbb Z^+$. In carrying out the   algorithm, one gets the feeling that 
the primes $\mathbb P$ are positioned in the sequence of integers precisely where they  need to be in order to generate the rest of $\mathbb Z^+$ by multiplication. Despite the sense of inevitability associated to the primes, many well-known theorems and conjectures related to prime numbers are very difficult (or seemingly impossible) to prove. 
As a way of sidestepping a direct attack on   prime number phenomena, mathematicians   often use a probabilistic approach based on the prime number theorem. 

We recall that if $\pi(x)$ denotes the number of prime numbers less than or equal to $x\in\mathbb R^+$, the classical {\it prime number theorem} states that as $x\to \infty$, 
\begin{equation}\label{pnt}
    \pi(x) \  \sim \  \frac{x}{\log x},
\end{equation}
where $\log x$ denotes the natural logarithm function. Thus the probability $\pi(n)/n$ that $n\in \mathbb Z^+$ is prime is asymptotically equal to $1/\log n$ as $n$ increases. A random model of prime distribution might assume the primes are modeled by a random sequence with probability $1/\log n$ of the $n$th term being prime, $n\geq 2$, or impose similar probabilistic measures on the prime gaps as in Cram\'{e}r's model  (see e.g. \cite{Tao, Tenenbaum}).

In recent work \cite{primemodel}, the  author of the present paper and their collaborators develop a {\it deterministic} model of prime number distribution based on facts and  observations about integer partitions. Inequalities that relate the partition norm \cite{SS-norm} and supernorm \cite{supernorm} statistics  impose constraints on the distribution of primes. Based on this, the authors of \cite{primemodel} deduce the $n$th prime gap may be modeled  by the   relation 
\begin{equation}\label{estimate}
    p_{n+1}-p_n\   =\  2\left\lceil \frac{d(n)}{2}\right\rceil,
\end{equation}
where $d(n)$ is the divisor function. While this partition model of prime gaps is formulated based on rather heavy-handed simplifying assumptions, it yields a surprisingly accurate estimate of prime   distribution. Below, we summarize the   model of primes proposed in     \cite{primemodel}. 

For $x\in \mathbb R$, let $\left\lceil x\right\rceil$ denote the  usual   {\it ceiling function}.  For $x\geq 0$, let $\left\lfloor x\right\rfloor$ denote the  {\it floor function}; however,  {\it we modify the usual floor function definition to let  $\left\lfloor x\right\rfloor:=0$ if $x<0$}.\\

\noindent {\bf Partition model of prime numbers formulated in \cite{primemodel}.} {\it The $n$th prime  number   $p_n,\  n\geq 1,$  is modeled by setting $p_1=2$ and for $n\geq 2$ by using the formula 
\begin{equation}\label{main}p_{n}\  =\  1\  +\  2\sum_{j=1}^{n-1}\left\lceil \frac{d(j)}{2}\right\rceil\  +\   \varepsilon(n),\end{equation}
where  $d(k)$ is the divisor function and $\varepsilon(k)$ is an explicit error   term  that is negligible by comparison. The case $\varepsilon(n):=0$ for all $n\geq 2$ is referred to in \cite{primemodel} as Model 1. The case 
\begin{equation}\label{model2}\varepsilon(n):=\left\lfloor \pi_2(p_{n-1})-2\gamma(n-1) \right\rfloor\end{equation}
is referred to  as Model 2, where $\pi_2(k)$ is the number of semiprimes less than or equal to $k\geq 1$, and $\gamma=0.5772...$ is the Euler--Mascheroni constant. For a third case, set  $p_1=2,\  p_2=3,$ and for $n\geq 3$ use an asymptotic approximation to the above error term,  
\begin{equation}\label{model2*} \varepsilon(n):=\left\lfloor (n-1)\cdot\left(\log \log(n-1)-2\gamma \right)\right\rfloor,\end{equation}
where $\log x$ denotes the natural logarithm, to yield a model variant referred to as Model 2*.}\\

Model 1 predicts both  the prime number theorem and  the twin prime conjecture \cite{primemodel}.  
In this study we experiment with  computational corrections to formulate a new model.\\

\noindent {\bf Partition model of prime numbers with   variable parameters (Model 3).} {\it Let $r,t\in\mathbb R$ with $r\geq 0,\  0\leq t\leq 1$. Let $M=0.2614...$ denote the Meissel-Mertens constant. Set $p_1=2,\  p_2=3$, and for $n\geq 3$ use equation \eqref{main} to model $p_n$ with the error term 
\begin{equation}\label{model3} \varepsilon(n)\  =\  \varepsilon(n;r,t) :=\  \left\lfloor (n+r-1)\log \log(n+r-1)-(n-1)(2\gamma-Mt) \right\rfloor\end{equation}
to  yield the   model variant we will refer to   as Model 3.}\\

Models 1, 2 and 2*, as well as Model 3 with $r=0$, produce the sequence of values $$2, 3, 5, 7, 11, 13, 17, 19, 23, 27,$$
and do not produce the    sequence of primes from there; but     almost produce the correct sequence of {\it prime gaps}.
See Table \ref{table1} to compare       estimates for   $\pi(n)$ from different models.

\vskip.1in

\begin{table}[h]    \centering
    \begin{tabular}{| c | c | c | c | c | c | c |}
\hline
 $n$\  & $\pi(n)$ & $n/\operatorname{log} n$  & Model 1 &  Model 2 &   Model 2*  & \begin{tabular}{@{}c@{}} Model 3 \\ $r=6, \  t=0.11$\end{tabular}  \\ \hline
 $10$\  & $4$ & $4.34...$   & 4  & 4 & 4 & 4\\ \hline
 $100$\  & $25$ & $21.71...$   & 27 & 26 & 27  & 25\\ \hline
  $1000$\  & $168$ & $144.76...$  & 184 & 168 & 171 & 168\\ \hline
   $10{,}000$\  & $1229$ & $1085.73...$   & 1352 & 1212 & 1233  & 1228\\ \hline
    $100{,}000$\  & $9592$ & $8685.88...$   & 10{,}602 & 9435 & 9618 & 9592 \\ \hline
    $1{,}000{,}000$\  & $78{,}498$ & $72{,}382.41...$   & 86{,}739 & 77{,}322  & 78{,}740 & 78{,}575\\ \hline
\end{tabular}
\vskip.1in
    \caption{Comparing   estimates for $\pi(n)$; Model 3 depends on  parameters $r, t$}
    \label{table1}
    \end{table}

\subsection
{Notations}
We define the notations required to  discuss these models. Let $\mathbb Z^+$ denote the   {\it positive integers}. Let  $\mathbb P$ denote the  {\it prime numbers}. 
Let $p_i\in\mathbb P$   denote the $i$th prime number, viz. $p_1=2, p_2=3, p_3=5,$ etc.; we define $p_0:=1$, and   refer to the subscript $i\in\mathbb Z^+$ of  $p_i\in\mathbb P$ as the {\it index} of the prime number.  
For $n\in \mathbb Z^+$,  we    write its {\it prime factorization} as $n=p_1^{a_1}p_2^{a_2}p_3^{a_3}\cdots p_i^{a_i}\cdots$,  $a_i\in \mathbb Z_{\geq 0}$,  where only finitely many primes $p_i$ have nonzero {\it multiplicity} (number of occurrences);   omit the factor $p_i$ from the notation if $a_i=0$. Define a {\it semiprime} to be a positive integer having exactly two prime factors, and define a {\it $k$-almost prime} to be  a positive integer with   $k\geq 2$ factors.  For $x\in \mathbb R^+$, let $\pi(x)$ denote the {\it prime counting function}, the number of primes less than or equal to $x$, and let $\pi_k(x)$ denote the {\it number of $k$-almost primes} less than or equal to $x$.

Let $\mathcal P$ denote the set of {\it integer partitions}, unordered finite multisets of positive integers including the empty partition $\emptyset\in \mathcal P$ (see e.g. \cite{And}).   
 For a nonempty partition $\lambda\in \mathcal P$, we write $\lambda=(\lambda_1, \lambda_2, \lambda_3, \dots, \lambda_r),\  \lambda_i\in\mathbb Z^+,$    $\lambda_1\geq \lambda_2\geq \dots \geq \lambda_r \geq 1$. 
For $\lambda \in \mathcal P$, let $|\lambda|\geq 0$ denote the {\it size} (sum of parts), let $\ell(\lambda)=r$ be the {\it length} (number of parts), 
and let $m_i=m_i(\lambda)\geq 0$ be the {\it multiplicity}   of $i\in \mathbb Z^+$ as a part of partition $\lambda$. Define a {\it non-unitary partition} to be a partition with no part equal to one.  

Let $\gamma=0.5772...$ denote the {\it Euler-Mascheroni constant} and let   $M=0.2614...$ denote the {\it Meissel-Mertens constant}. 

A  new, multiplicative partition statistic is introduced in \cite{Robert_zeta}.  
The {\it norm} $N(\lambda)$ of   partition $\lambda = (\lambda_1, \lambda_2, \lambda_3, \dots, \lambda_r)$ is the product of its parts:
\begin{equation}\label{normdef}
N(\lambda)\  :=\  \lambda_1  \lambda_2  \lambda_3 \cdots  \lambda_r\  =\  1^{m_1} 2^{m_2} 3^{m_3} \cdots  i^{m_i}\cdots \in \mathbb Z^+.
\end{equation}
We define $N(\emptyset):=1$ (it is an empty product); see \cite{SS-norm} for more on the partition norm. 

Furthermore, in \cite{supernorm} another multiplicative partition statistic is defined. The {\it supernorm} $\widehat{{N}}$ of  partition $\lambda = (\lambda_1, \lambda_2, \lambda_3, \dots, \lambda_r)$ is the product 
 \begin{equation}\label{supernormdef}
 \widehat{N}(\lambda)\  :=\  p_{\lambda_1}p_{\lambda_2}\cdots p_{\lambda_r} \  =\  2^{m_1} 3^{m_2} 5^{m_3} \cdots p_i^{m_i}\cdots \in \mathbb Z^+,\end{equation}
where $p_i\in \mathbb P$ is the $i$th prime number, $i\geq 1$, and $m_j=m_j(\lambda)\geq 0$. We define $\widehat{N}(\emptyset):=1$; see \cite{Lagarias, Lagarias2} for further reading about the partition supernorm. 

\section{Improving the computational model of  $\pi(n)$}\label{sect2}

\subsection{Brief summary of the ideas behind Models 1, 2 and 2*}\label{model}
In \cite{primemodel}, Model 1 follows from two elementary observations and a subsequent simplifying assumption. \\ 

\noindent {\bf Observation 1.} {\it The supernorm induces a bijection between  partitions and the prime factorizations of positive integers. The integers in prime gaps are  the images of certain subsets of  $\mathcal P$ under the supernorm map.  Thus, as is noted in \cite[Prop. 2.1]{primemodel}, {measuring the $n$th prime gap is equivalent to enumerating     partitions that map into the interval $[p_n, p_{n+1})$  under the supernorm}:
\begin{equation}\label{gapequiv}
p_{n+1}-p_n \  \  =\  \  \    \#\left\{\lambda\in\mathcal P\  :\  p_n \leq \widehat{N}(\lambda) <p_{n+1}\right\}.\end{equation} }

\noindent 
{\bf Observation 2.} {\it Based on a mix of proof and computational evidence, it is observed in \cite[Eqs. 1.9--1.10]{primemodel} that if $\lambda$ is a {non-unitary partition}   and $N(\lambda)\geq 5$,  
\begin{equation}\label{normestimatethm4} N(\lambda)\  <\  p_{N(\lambda)}\  \leq\  \widehat{N}(\lambda)\  \leq \  N(\lambda)^{{\log 3}/{\log 2}}.\end{equation}
Noting $\frac{\log 3}{\log 2} = 1.5849...$, the supernorm maps $\lambda$ to a   ``small'' interval immediately above the prime $p_{N(\lambda)}$. Moreover, the fewer parts $\lambda$ has, the closer we expect $\widehat{N}(\lambda)$ is to $p_{N(\lambda)}$.}\\

Observation 1 suggests a reasonable way to estimate the size $p_{n+1}-p_n$  of the $n$th prime gap is to guess which partitions map into the interval under $\widehat{N}$, and enumerate them. Observation 2 suggests a good place to start might be with partitions whose images under $\widehat{N}$ we already expect to be just above $p_n$, viz. non-unitary partitions of norm $n$ with  few parts. These partitions map to odd integers we anticipate to lie  just above $p_n$ on the number line. In \cite{primemodel}, a somewhat drastic simplification  is assumed based on this reasoning. \\

\noindent {\bf Simplifying Assumption.} 
{\it   Following \cite{primemodel}, assume for simplicity that {all} odd integers in the interval $[p_n, p_{n+1}),\  n\geq 2$,  are the images of non-unitary partitions with norms equal to $n$,   having    one or two parts, under the supernorm $\widehat{N}$. That is, assume for simplicity  
\begin{equation}\label{gapequiv2}
p_{n+1}-p_n \  \  =\  \  2\cdot  \#\left\{\lambda\in \mathcal P\  :\   \  N(\lambda)=n, \  m_1(\lambda)=0,\    \ell(\lambda)=1\  \text{or}\  2\right\};\end{equation}
the factor of $2$ on the right   doubles the number of odd primes and semiprimes in the interval to count the   even integers that follow them.}\\ 

As   noted in \cite{primemodel}, assumption   \eqref{gapequiv2} is equivalent   to writing    
\begin{equation}\label{estimate2}
    p_{n+1}-p_n\  =\  2\left\lceil \frac{d(n)}{2}\right\rceil.
\end{equation}
Each length-two partition $\lambda=(\lambda_1, \lambda_2)$  such that $N(\lambda)=\lambda_1 \lambda_2=n$ represents a pair of divisors $\lambda_1$ and $\lambda_2$ of $n$, and we associate the one-part partition $(n)$   with the divisor pair $n$ and $1$, justifying \eqref{estimate} above.  
Then for $n\geq 2$, one obtains from \eqref{estimate2} by telescoping series 
\begin{equation}
    p_n \ =\  p_2+\sum_{k=2}^{n-1}(p_{k+1}-p_{k})\  =\  3+2\sum_{j=2}^{n-1}\left\lceil \frac{d(j)}{2}\right\rceil\  =\  1+2\sum_{j=1}^{n-1}\left\lceil \frac{d(j)}{2}\right\rceil,
\end{equation}
since $3=1+2\left\lceil {d(1)}/{2}\right\rceil$. This is the $\varepsilon(n)=0$ case of \eqref{main} that defines Model 1 in \cite{primemodel}.

\begin{remark}
    Note that   $2\left\lceil d(k)/{2}\right\rceil=d(k)$ for all $k\neq m^2, m\in \mathbb Z^+$. If $k$ is a perfect square then $2\left\lceil {d(k)}/{2}\right\rceil=d(k)+1$. Then for $n\geq 2$, one may simplify computations by using 
\begin{equation}\label{simplify}
2\sum_{k=1}^{n-1}\left\lceil \frac{d(k)}{2}\right\rceil\  =\   \sum_{k=1}^{n-1} d(k)\  +\   \left\lfloor \sqrt{n-1}\right\rfloor\  =\  \sum_{k=1}^{n-1}\left\lfloor \frac{n-1}{k}\right\rfloor\  +\    \left\lfloor \sqrt{n-1}\right\rfloor.
\end{equation}
\end{remark}

Since integers in prime gaps with more than two prime factors are not enumerated in Model 1, it is anticipated that the model should underestimate the prime gap sizes, and thus overestimate $\pi(n)$; Table \ref{table1} suggests this is indeed the case. 

In \cite{primemodel},  Model 2 is derived from Model 1 by making a simplistic estimate of the missing contribution  to prime gaps of partitions having {\it three} parts and  norms equal to $n$.  
The $\pi_2(p_{n-1})$ term in \eqref{model2} arises  in this estimate; the   $2\gamma(n-1)$  term is subtracted for computational reasons and    significantly improves the model's estimates. See \cite{primemodel} for details.

Model 2* is then produced in \cite{primemodel} by substituting the following   estimate, derived from a formula of Landau \cite{handbuch}, into Model 2 to yield the estimated  correction term   $\varepsilon(n)$   in \eqref{model2*}: 
\begin{equation}\label{pi_2}
\pi_2(p_{n-1})\sim (n-1)\log\log (n-1)\end{equation} 
as $n\to \infty$. The authors of  \cite{primemodel} make a     comparison of Models 2 and 2*: ``Since it involves simpler computations and is  evidently    more accurate in its estimates of $\pi(n)$ as $n$ increases, we consider Model 2* to be a better model of primes.''

In the following section,   we focus on   
refining the error term \eqref{pi_2}.

\subsection{Refining the estimate of $\pi_2(p_{n-1})$ in the error term (Model 3)} 
In \cite{primemodel}, the authors note that as $n$ increases, ``the missing contribution of partitions (and   factorizations) of lengths greater than three   creeps up, causing the [Model 2 and 2*] estimates to become farther from exact as $n$ increases.'' It is reasonable to expect that estimating the contributions of partitions with more than three parts may improve the models' accuracy; Model 2 was formulated for the same reason, but related to partitions of lengths greater than two. That the Model 2 correction term \eqref{model2} and subsequent Model 2* correction \eqref{model2*} were suggested by the premises of Model 1, and {\it do} produce improved estimates of prime distribution, strengthens the plausibility of the model's premises to some extent.  

However, as the number of parts increases, the   premise is compromised that partitions of norm $n$ with {\it few parts}  should map immediately above $p_n$. 
To eliminate this concern, let us focus instead on modifying the asymptotic estimate \eqref{pi_2} to produce a better fit.

\subsubsection{Interpolating between Models 2 and 2*}
As discussed in \cite{primemodel},  inspection of Table \ref{table1}  suggests that from some point on, the Model 2 estimate for $\pi(n)$ is strictly less than $\pi(n)$ itself while the Model 2* estimate is strictly greater than $\pi(n)$; this observation is recorded in  \cite[Conjecture 4.1]{primemodel}. {\it If these inequalities were true, it would imply that Model 2 accurately captures the behavior of $\pi(n)$ with  smaller error than the error in the Landau estimate \eqref{pi_2}
 for $\pi_2(p_{n-1})$.} (In \cite{primemodel}, the authors note this   conjecture relies on limited data.)

It is proved in \cite{primemodel} based on a formula from  \cite{MM} that, for $M$   the Meissel-Mertens constant,  
\begin{equation}\label{underestimate}
    \pi_2(p_{n})-n\log\log n\  \sim\    M n
\end{equation} 
as $n\to \infty$. Then   a better approximation to \eqref{model2} than $\eqref{model2*}$ provides  is given  by 
\begin{equation}\label{model2.1} \varepsilon(n):=\left\lfloor (n-1)\cdot\left(\log \log(n-1)-2\gamma+M \right)\right\rfloor.\end{equation}
One can effectively interpolate between Model 2 and Model 2* by introducing a parameter $t\in [0, 1]$ into \eqref{model2.1}, so that $t=0$ gives Model 2* and $t=1$   approximates  Model 2:
\begin{equation}\label{model3.1} \varepsilon(n)=\varepsilon(n;t):=\left\lfloor (n-1)\cdot\left(\log \log(n-1)-2\gamma+Mt \right)\right\rfloor.\end{equation}
By varying $t$, perhaps one can locate a ``sweet spot''   between the estimates of Model 2 and Model 2* to produce a better model of $\pi(n)$. See Table \ref{table2} for examples.

\subsubsection{Introducing a second parameter to produce a better fit}
Another variable parameter can be introduced into the asymptotic estimate $\pi_2(p_{n})\sim n\log\log n$ in a somewhat {\it ad hoc} manner, to    further fine tune estimates for $\pi(n)$. In \cite{EliRobert} it is noted that for   fixed $r\geq 0$,  it is also the case that $(n+r)\log\log (n+r)\sim n \log \log n$ as $n\to \infty$.\footnote{The author is    grateful to Eli DeWitt who performed    computational experiments   as an undergraduate research student at Michigan Technological University (see \cite{EliRobert}) that informed this study.}  
Computational experiments reported in \cite{EliRobert}   show that   replacing \eqref{model2*} in Model 2* with 
\begin{equation}\label{model2.3} 
\varepsilon(n) = \varepsilon(n;r):=\left\lfloor (n+r-1)\log \log(n+r-1)-2\gamma(n-1)\right\rfloor
\end{equation}
and varying $r$ may improve Model 2* estimates of $\pi(n)$ at small $n$, with a decreasing effect as $n\to\infty$ since  $(n+r)\log\log (n+r)- n \log \log n\sim 0$. See Table \ref{table3} for examples.

\begin{table}[h]    \centering
    \begin{tabular}{| c | c | c | c | c | c | }
\hline
 $n$\  & $\pi(n)$ & $t=0$  & $t=0.1$ &  $t=0.5$ &  $t=1$   \\ \hline
 $10$\  & $4$ & $4$   & 4  & 4 & 4  \\ \hline
 $100$\  & $25$ & $27$   & 27 & 26 & 26   \\ \hline
  $1000$\  & $168$ & $171$  & 170 & 168 & 165 \\ \hline
   $10{,}000$\  & $1229$ & $1233$   & 1230 & 1216 & 1200  \\ \hline
    $100{,}000$\  & $9592$ & $9618$   & 9595 & 9510 & 9404   \\ \hline
    $1{,}000{,}000$\  & $78{,}498$ & $78{,}740$   & $78{,}591$ & 78{,}000  & 77{,}274  \\ \hline
\end{tabular}
\vskip.1in
    \caption{Estimates for $\pi(n)$ when  varying $t$ in equation \eqref{model3.1}}
    \label{table2}
    \end{table}

\subsubsection{Combining the above modifications} 
We    combine the modifications introduced to $\varepsilon(n)$ in \eqref{model3.1}
 and \eqref{model2.3}
  to yield the new   model we mentioned in Section \ref{sect1}, as follows.\\

\noindent {\bf Model 3.} {\it The prime numbers $p_1, p_2, p_3,$ $\dots, $ can be modeled by the sequence having initial values $p_1=2,\  p_2=3$,   and for $n\geq 3$  having the values  
$$p_{n}\  =\  1 \  +\   2\sum_{k=1}^{n-1}\left\lceil \frac{d(k)}{2}\right\rceil\  +\  \left\lfloor (n+r-1)\log \log(n+r-1)-(n-1)(2\gamma-Mt) \right\rfloor,$$}
with $r,t\in\mathbb R,\  r\geq 0,\  0\leq t\leq 1$.\\

See Table \ref{table1} to compare the Model 3 estimates for $\pi(n)$ to estimates from   other models. 

Setting $r=t=0$, Model 3 reduces to Model 2* from \cite{primemodel}. If $r=0$, then computations show Model 3 outputs the same initial values  $2, 3, 5, 7, 11, 13, 17, 19,  23, 27,$ as the previous iterations of the model. However, as $r>0$ increases, the   terms of the modeled $p_n$ sequence become spaced farther  apart and Model 3 does not model the sequence of primes at all, while still providing a reasonably good model of  $\pi(n)$ that  remains asymptotic to $n/\log n$. 

Here is an example of this phenomenon for a particular selection of the values $r,t$. Note   in Table \ref{table2} that the $t=0.1$ column provides quite a good fit to $\pi(n)$; also,  setting $r=6$   is noted to yield a good fit in \cite{EliRobert}. One finds   that setting   $r=6,\  t=0.11,$  in Model 3   produces   an accurate estimate of $\pi(n)$   up to at least $n=1{,}000{,}000$, as Table \ref{table1} displays. But setting $r=6, \  t=0.11,$ in Model 3 outputs the sequence of values $2, 3, 8, 10, 14, 16, 21, 23, ...$. This sequence has the wrong values from the third term on, including   even integers that are unmistakably composite,\footnote{In \cite{primemodel} it is noted that replacing the floor function $\left\lfloor x\right\rfloor$ in the correction terms of these models with a variant $\left\lfloor x\right\rfloor^*$ such that $\left\lfloor x\right\rfloor^*=\left\lfloor x\right\rfloor$ if $\left\lfloor x\right\rfloor$ is even, and $\left\lfloor x\right\rfloor^*=\left\lfloor x\right\rfloor-1$ if $\left\lfloor x\right\rfloor$ is odd, will produce a sequence of only odd numbers without affecting estimates for $\pi(n)$.} yet is almost equinumerous  with the primes up to one million.

\section{Closing remarks}
The    models  formulated in \cite{primemodel} and here, as well as   recent works like \cite{Ono1, supernorm}, demonstrate that   partition theory -- the theory of addition -- encodes     information about prime numbers. To what extent this insight is useful to the theory of prime distribution, compared to classical approaches like random models, is yet to be seen. In \cite{Ono1},  prime numbers are proved to be intrinsically connected to the coefficients of certain quasimodular forms via special partition-theoretic functions; this application is seemingly unconnected to the topic of prime distribution. 
In \cite{primemodel}, while a partition model of primes and of the prime number theorem is established, the authors test new predictions from Model 1 with mixed findings and, in   successful tests, with results   only weakly supporting predictions. Perhaps the prime number formulas produced by these models are   as accurate as they are due to entirely   different reasons  from those given in \cite{primemodel}. 
Yet because these partition   models for $\pi(n)$ are  accurate at small   $n$ and    provide   discrete models of the prime number theorem as $n$ increases, one wonders:  do they     capture something essential about prime distribution?

\vskip.1in

\begin{table}[h]    \centering
    \begin{tabular}{| c | c | c | c | c | c | c |}
\hline
 $n$\  & $\pi(n)$   & $r=0$  & $r=1$ &  $r=6$ &  $r=10$  & $r=20$ \\ \hline
 $10$\  & $4$ & $4$   & 4  & 4 & 2 & 2\\ \hline
 $100$\  & $25$ & $27$   & 27 & 25 & 24  & 20\\ \hline
  $1000$\  & $168$ & $171$  & 170 & 169 & 168 & 166\\ \hline
   $10{,}000$\  & 1229& $1233$   & 1232 & 1232 & 1232  & 1229\\ \hline
    $100{,}000$\  & 9592 & $9618$   & 9618 & 9616 & 9615& 9614 \\ \hline
    $1{,}000{,}000$\  & $78{,}498$ & $78{,}740$   & 78{,}740 & 78{,}740  & 78{,}738 & 78{,}736\\ \hline
\end{tabular} \vskip.1in
    \caption{Estimates for $\pi(n)$  when  varying $r$ in \eqref{model2.3},   reproduced from \cite{EliRobert}}
    \label{table3}
    \end{table}

\section*{Acknowledgments}
For calculations, the author used   Python programming language. I am thankful  to Eli DeWitt for performing computational tests  of Model 2* that informed the formulation of Model 3 during a summer REU  supported by the Department of Mathematical Sciences at Michigan Technological University, and  to Aidan Botkin for sharing computer code with me that was used to produce Models 2 and 2* in \cite{primemodel}.    

\end{document}